\documentclass[reqno,11pt]{amsart}
\usepackage{amssymb}
\usepackage{amsmath}
\usepackage{times}
\usepackage{graphicx}
\usepackage{hyperref}

\newtheorem{theorem}{Theorem}[section]
\newtheorem{lemma}[theorem]{Lemma}
\newtheorem{proposition}[theorem]{Proposition}
\newtheorem{definition}[theorem]{Definition}
\newtheorem{corollary}[theorem]{Corollary}

\begin{document}
\title{On Pfaffian random point fields}
\thanks{Statistical Laboratory, Department of Mathematics, University of
Cambridge, Cambridge, UK. e-mail: vladislav.kargin@gmail.com}
\date{May 2013}
\author{V. Kargin}
\maketitle

\begin{center}
\textbf{Abstract}
\end{center}

\begin{quotation}
We study Pfaffian random point fields by using the Moore-Dyson quaternion
determinants. First, we give sufficient conditions that ensure that a
self-dual quaternion kernel defines a valid random point field, and then we
prove a CLT for Pfaffian point fields. The proofs are based on a new
quaternion extension of the Cauchy-Binet determinantal identity. In addition, we derive the Fredholm determinantal formulas for the Pfaffian point fields which use the quaternion determinant.
\end{quotation}

\bigskip

\section{Introduction}

A determinantal random point field is a random collection of points with the
probability distribution that can be written as a determinant. The
determinantal point fields describe various mathematical objects including
eigenvalues of random matrices, zeros of random analytic functions,
non-intersecting random paths, and spanning trees on networks. It is
conjectured that they are also related to other important objects such as
Riemann's zeta function zeros and the spectrum of chaotic dynamical systems.
See \cite{soshnikov00b}, \cite{hkpv06}, \cite{hkpv09} and \cite{borodin09}
for reviews.

The definition of the determinantal point field uses the standard matrix
determinant. However, in some applications, the distribution of a random
point field can be represented as a determinant of a quaternion matrix. An
important example is provided by eigenvalues of orthogonal and symplectic
random matrix ensembles. It is natural to study these fields as a
generalization of the usual determinantal point fields and study their
properties.

While this generalization can be thought of as a quaternion determinantal
point field, it is also equivalent to the \emph{Pfaffian} \emph{random point
field}, which was recently studied in \cite{rains00}, \cite{soshnikov03},
and \cite{borodin_rains05}. In order to give a precise definition, we recall
some preliminary notations.

A \emph{random point field} $\mathcal{X=}\left( X,\mathcal{B},\mathbb{P}%
\right) $ on a measurable space $\Lambda $ is a probability measure $\mathbb{%
P}$ on the space $X$ of all possible countable configurations of points in $%
\Lambda $. It is convenient to think about this object as a collection of
functions that sends every $n$-tuple of non-negative integers $\left(
k_{1},\ldots ,k_{n}\right) $ and every $n$-tuple of measurable subsets of $%
\Lambda ,$ $\left( A_{1},\ldots ,A_{n}\right) $ to a non-negative number,
which can be interpreted as a probability to find $k_{i}$ points in the set $%
A_{i}.$ These functions are to satisfy some consistency conditions, which we
do not specify here. The reader is advised to consult paper \cite{lenard73}
or book \cite{daley_vere-jones88} for more detail. We will also call the
process $\mathcal{X}$ a $\Lambda $\emph{-valued random point field}.

Suppose that $\Lambda $ is a space with measure $\mu .$ Let $\#\left(
A_{i}\right) $ denote the number of points of $\mathcal{X}$ located in the
set $A_{i}.$

\begin{definition}
A locally integrable function $R_{k}$: $\Lambda ^{k}\rightarrow \mathbb{R}%
_{+}^{1}$ is called a $k$\emph{-point correlation function} of a random
point field $\mathcal{X=}\left( X,\mathcal{B},\mathbb{P}\right) $ with
respect to the measure $\mu ,$ if for any disjoint measurable subsets $%
A_{1},\ldots ,A_{m}$ of $\Lambda $ and any non-negative integers $%
k_{1},\ldots ,k_{m},$ such that $\sum_{i=1}^{m}k_{i}=k,$ the following
formula holds: 
\begin{equation}
\mathbb{E}\prod_{i=1}^{m}\left[ \#\left( A_{i}\right) \ldots \left( \#\left(
A_{i}\right) -k_{i}+1\right) \right] =\int_{A_{1}^{k_{1}}\times \ldots
\times A_{m}^{k_{m}}}R_{k}\left( x_{1},\ldots ,x_{k}\right) d\mu \left(
x_{1}\right) \cdots d\mu \left( x_{k}\right) ,
\label{def_correlation_formulas}
\end{equation}%
where $\mathbb{E}$ denote expectation with respect to measure $\mathbb{P}$.
\end{definition}

Note that on the left is the expected number of ordered configurations of
points such that set $A_{i}$ contains $k_{i}$ points. Note also that the
correlation functions are defined only up to sets of measure $0.$

\begin{definition}
\label{definition_main}A $\Lambda $-valued random point field $\mathcal{X}$
is called a $\emph{Pfaffian}$\emph{\ random point field} if its correlation
functions can be written as quaternion determinants:%
\begin{equation}
R_{m}\left( x_{1},\ldots ,x_{m}\right) =\mathrm{Det}_{M}\left. \left(
K\left( x_{i},x_{j}\right) \right) \right\vert _{1\leq i,j\leq m},\text{ }%
x_{1},\ldots ,x_{m}\in \Lambda ,\text{ }m=1,2,...,
\label{def_correlation_determinant}
\end{equation}%
where $K\left( x,y\right) $ is a self-dual quaternion kernel (that is, $%
K\left( y,x\right) =\left( K\left( x,y\right) \right) ^{\ast }$), and $%
\mathrm{Det}_{M}$ is the Moore-Dyson quaternion determinant.
\end{definition}

In this definition the function $K\left( x,y\right) $ takes value in the
algebra of complexified quaternions $\mathbb{Q}_{\mathbb{C}}.$

(Recall that real quaternions can be written $q=s+x\mathbf{i}+y\mathbf{j}+z%
\mathbf{k,}$ where $s,x,y,$ and $z$ are real and where $\mathbf{i},\mathbf{j}%
,\mathbf{k}$ denote the quaternion units with the rules $\mathbf{ij}=-%
\mathbf{ji}=\mathbf{k}$\ and so on. The complexified quaternions are allowed
to have complex coefficients $s,x,y,$ and $z.$ For both real and
complexified quaternions, the conjugate of $q$ is defined as $q^{\ast }=s-x%
\mathbf{i}-y\mathbf{j}-z\mathbf{k}.$ In this paper when we say quaternions,
we mean complexified quaternions, and we say real quaternions for
quaternions with real coefficients. Quaternion matrices are matrices whose
entries are quaternions. The \emph{dual} of a quaternion matrix $X$ is
defined as a matrix $X^{\ast },$ for which $\left( X^{\ast }\right)
_{lk}=\left( X_{kl}\right) ^{\ast }.$ \emph{Self-dual} quaternion matrices
are defined by the property that $X^{\ast }=X.$)

The name Pfaffian comes from a different definition of this object given in 
\cite{soshnikov03} and \cite{borodin_rains05}. \ They call a random point
field Pfaffian if there exists a $2\times 2$ matrix-valued skew-symmetric
kernel $K$ on $X$ such that the correlation functions of the process have
the form

\begin{equation*}
R_{m}\left( x_{1},\ldots ,x_{m}\right) =\mathrm{Pf}\left[ K\left(
x_{i},x_{j}\right) \right] _{i,j=1}^{m},\text{ }x_{1},\ldots ,x_{m}\in X,%
\text{ }m=1,2,\ldots
\end{equation*}%
(The notation $\mathrm{Pf}$ in the right-hand side stands for the Pfaffian.)
This definition is equivalent to ours. Indeed, one can take the complex
matrix representation of the quaternion kernel $K\left( x_{i},x_{j}\right) $
and write the quaternion determinant in Definition \ref{definition_main} as
a Pfaffian of a 2-by-2 matrix-valued skew-symmetric kernel equal to this
representation multiplied by a matrix $J=\left( 
\begin{array}{cc}
0 & I \\ 
-I & 0%
\end{array}%
\right) $. (See formula (\ref{formula_determinant_pfaffian}) in Appendix.)

In the opposite way, suppose that the field is Pfaffian with $2\times 2$
matrix-valued kernel $K\left( x,y\right) $. Since 
\begin{equation*}
K\left( x,y\right) =\left( 
\begin{array}{cc}
A\left( x,y\right) & B\left( x,y\right) \\ 
-B\left( y,x\right) & D\left( x,y\right)%
\end{array}%
\right)
\end{equation*}%
is skew-symmetric, hence $A\left( y,x\right) =-A\left( x,y\right) $ and $%
D\left( y,x\right) =-D\left( x,y\right) $ and it is easy to check that 
\begin{equation*}
\left( 
\begin{array}{cc}
B\left( y,x\right) & -D\left( x,y\right) \\ 
A\left( x,y\right) & B\left( x,y\right)%
\end{array}%
\right)
\end{equation*}%
is a complex representation of a self-dual quaternion kernel $\widetilde{K}%
\left( x,y\right) .$ By using the relation between the quaternion
determinant and the Pfaffian, we conclude that the correlations can be
written as quaternion determinants of the kernel $\widetilde{K}\left(
x,y\right) .$

Here are some examples of Pfaffian fields.

\textbf{Example 1.} It was shown by Dyson \cite{dyson70} that the point
field of eigenvalues of the circular orthogonal and circular symplectic
ensembles of random matrices have correlation functions that can be written
as quaternion determinants. Hence they form a Pfaffian point field that
takes values in the unit circle of the complex plane. Later this result was
extended by Mehta (see Chapters 7 and 8 in \cite{mehta04}) to the case of
Gaussian orthogonal and symplectic ensembles. In this case, the eigenvalues
form a real-valued Pfaffian point field. We will say more about the
symplectic random matrix ensembles in the last section.

\textbf{Example 2.} (Ginibre on $\mathbb{Q}$) Let $\Lambda $=$\mathbb{Q}$,
where $\mathbb{Q}\simeq \mathbb{R}^{4}$ denotes real quaternions. Take the
background measure $d\mu \left( z\right) =\pi ^{-2}e^{-\left\vert
z\right\vert ^{2}}dm\left( z\right) ,$ where $dm\left( z\right) $ is the
Lebesgue measure on $\mathbb{Q}$, and define the kernel 
\begin{equation}
K_{n}\left( z,w\right) =\sum_{k=0}^{n}\frac{z^{k}\left( w^{\ast }\right) ^{k}%
}{\left( k+1\right) !},\text{ where }z,w\in \mathbb{Q}
\label{quaternion_Ginibre}
\end{equation}%
This kernel corresponds to a Pfaffian point field that takes value in real
quaternions. The fact that kernel (\ref{quaternion_Ginibre}) defines a valid
random point field follows from Proposition \ref{prop_existence_Bernoulli}
below.

Recall for comparison that the usual Ginibre random point field is a
complex-valued point field of eigenvalues of a random $n$-by-$n$ complex
Gaussian matrix. It is determinantal with the kernel 
\begin{equation*}
K_{n}\left( z,w\right) =\sum_{k=0}^{n}\frac{z^{k}\left( \overline{w}\right)
^{k}}{k!},
\end{equation*}%
where $z$ and $w$ are in $\mathbb{C}$ and the background measure is $d\mu
\left( z\right) =\pi ^{-1}e^{-\left\vert z\right\vert ^{2}}dm\left( z\right)
.$ (See \cite{ginibre65}, Section 15.1 in \cite{mehta04}, and Section 4.3.7
in \cite{hkpv09}.)

\textbf{Example 3.} (Pfaffian Ginibre on $\mathbb{C}$) \ This is another
generalization of the Ginibre random point field. See Section 15.2 in \cite%
{mehta04} for details. It is a complex-valued Pfaffian point field. Let 
\begin{equation*}
\phi _{N}\left( u,v\right) =\frac{1}{2\pi }\sum_{0\leq i\leq j\leq N-1}\frac{%
2^{j}j!}{2^{i}i!}\frac{1}{\left( 2j+1\right) !}\left(
u^{2i}v^{2j+1}-v^{2i}u^{2j+1}\right) ,
\end{equation*}%
and define the quaternion kernel by its complex matrix representation: 
\begin{equation*}
\varphi \left( K_{N}\left( z,w\right) \right) =\left( 
\begin{array}{cc}
\phi _{N}\left( w,\overline{z}\right) & \phi _{N}\left( \overline{w},%
\overline{z}\right) \\ 
\phi _{N}\left( z,w\right) & \phi _{N}\left( z,\overline{w}\right)%
\end{array}%
\right) .
\end{equation*}%
(The map $\varphi $: $\mathbb{Q}_{\mathbb{C}}\rightarrow M_{2}\left( \mathbb{%
C}\right) $ is a bijection between complexified quaternions and the $2$-by-$%
2 $ complex matrices$.$ Its definition is standard and given in Appendix.)
Then this kernel defines a Pfaffian point field with respect to the signed
background measure $d\mu \left( z\right) =e^{-\left\vert z\right\vert
^{2}}\left( z-\overline{z}\right) dm\left( z\right) .$

\textbf{Example 4} (Bergman kernel on $\mathbb{Q}$). Let $\Lambda $ be the
unit disc in $\mathbb{Q}$ with the background measure $d\mu \left( z\right)
=\pi ^{-2}dm\left( z\right) ,$ where $dm\left( z\right) $ is the Lebesgue
measure on $\mathbb{Q\simeq R}^{4}$. Define the kernel 
\begin{equation}
K_{n}\left( z,w\right) =\sum_{k=0}^{n}\left( k+2\right) z^{k}\left( w^{\ast
}\right) ^{k},\text{ where }z,w\in \mathbb{Q}\text{.}
\end{equation}%
Then this kernel corresponds to a Pfaffian point field that takes value in
the unit disc of real quaternions.

For comparison, the Bergman kernel on the unit disc in $\mathbb{C}$ is given
by 
\begin{equation*}
K_{n}\left( z,w\right) =\sum_{k=0}^{n}\left( k+1\right) z^{k}\left( 
\overline{w}\right) ^{k},
\end{equation*}%
and corresponds to the determinantal point field of zeros of power series
with i.i.d complex Gaussian coefficients. (See Section 15.2 in \cite{hkpv09}%
).

Which self-dual quaternion kernels correspond to random point fields?

Note that a correlation function is automatically symmetric if it is defined
as in (\ref{def_correlation_determinant}). That is, we have 
\begin{equation*}
R_{m}\left( x_{\sigma (1)},\ldots ,x_{\sigma (m)}\right) =R_{m}\left(
x_{1},\ldots ,x_{m}\right)
\end{equation*}%
for every permutation $\sigma \in S_{m}.$ Hence, by the Lenard criterion (%
\cite{lenard75a} and \cite{lenard75b}), positivity is a necessary and
sufficient condition that ensures that a kernel corresponds to a random
point field. This condition can be explained as follows. Let $X$ be the
space of configurations of points in $\Lambda .$ Let $\varphi =\{\varphi
_{k}\}$ denote an arbitrary sequence of real-valued functions $\varphi _{k}$
over $\Lambda ^{k}$ which have the property that they are zero if at least
one of its arguments is outside of a certain compact set in $\Lambda .$
Define operator $S$ on $\varphi $ as follows: $S$ maps sequence $\varphi $
to a function $S\varphi $ on $X$ 
\begin{equation*}
\left( S\varphi \right) \left( x\right) =\sum \varphi _{k}\left(
x_{i_{1}},\ldots ,x_{i_{k}}\right)
\end{equation*}%
where $x=(x_{1},x_{2},\ldots )$ is an enumeration of a point configuration
in $X,$ and the sum is extended over all finite sequences of distinct
positive integers $\left( i_{1},\ldots ,i_{k}\right) $ including the empty
sequence. (This sum is essentially finite since the definition of the space
of configurations requires that the number of points of any configuration in
any compact subset of $\Lambda $ be finite.) The positivity condition says
that if $\left( S\varphi \right) \left( x\right) \geq 0$ for all $x\in X$
(including the empty sequence)$,$ then it must be true that 
\begin{equation*}
\int \varphi \left( x\right) d\rho :=\varphi
_{0}+\sum_{k=1}^{N}\int_{\Lambda ^{k}}\varphi _{k}\left( x_{1},\ldots
,x_{k}\right) R_{k}\left( x_{1},\ldots ,x_{k}\right) d\mu \left(
x_{1}\right) \ldots d\mu \left( x_{k}\right) \geq 0.
\end{equation*}%
This criterion is useful for fields with small number of points. However, in
general this criterion is often difficult to verify.

In the case of the usual determinantal fields with self-adjoint kernels a
much simpler criterion is given by the Macchi-Soshnikov theorem (see \cite%
{macchi75} and \cite{soshnikov00b}) that says that a self-adjoint kernel $%
K\left( x,y\right) $ defines a determinantal point field if and only if the
corresponding operator $K$ is in the trace class and all its eigenvalues are
in the interval $\left[ 0,1\right] .$

Our first result is a weaker version of this theorem for Pfaffian point
fields.

Suppose that a quaternion kernel $K$ can be written as follows: 
\begin{equation}
K\left( x,y\right) =\sum_{k=1}^{\infty }\lambda _{k}u_{k}\left( x\right)
u_{k}^{\ast }\left( y\right) ,  \label{def_quaternion_kernel}
\end{equation}%
where $\lambda _{k}$ are scalar and $u_{k}\left( x\right) $ is an
orthonormal system of quaternion functions: 
\begin{equation*}
\int_{\Lambda }u_{k}^{\ast }\left( x\right) u_{l}\left( x\right) d\mu
(x)=\delta _{kl}.
\end{equation*}%
(The series in (\ref{def_quaternion_kernel}) are assumed to be absolutely
convergent almost everywhere.) We will say in this case that $K\left(
x,y\right) $ has a \emph{diagonal form} with eigenvalues $\lambda _{k}.$ If $%
\lambda _{k}$ are real and $u_{k}\left( x\right) $ is an orthonormal system
of \emph{real quaternion} functions, then we will say that $K\left(
x,y\right) $ has a \emph{real diagonal form}$.$

\begin{theorem}
\label{theorem_existence}Suppose that a quaternion kernel $K\left(
x,y\right) $ has a real diagonal form with eigenvalues $\lambda _{k}.$
Assume that all $\lambda _{k}\in \left[ 0,1\right] ,$ and that $%
\sum_{k=1}^{\infty }\lambda _{k}<\infty .$ Then the kernel $K\left(
x,y\right) $ defines a Pfaffian point field with finite expected number of
points.
\end{theorem}

The conditions of this theorem are sufficient but not necessary since there
exist self-dual quaternion kernels which do not have a real diagonal form
and still define a valid point field. For example, let the space $\Lambda $
consist of two points and has the counting background measure. Let $K=\frac{1%
}{2}\left( 
\begin{array}{cc}
1 & -a \\ 
a & 1%
\end{array}%
\right) ,$where $a=\left( 3i\mathbf{i-}5\mathbf{j}\right) /4$ so that $%
a^{2}=-1.$ This matrix is self-dual with determinant zero. It defines a
random point field that has exactly one point uniformly distributed on $%
\Lambda .$ This kernel has a diagonal form but it does not have a real
diagonal form. (The eigenvector is a complex quaternion vector.)

Here is another example. Let $\Lambda $ consist of two points, and let $%
K=\left( 
\begin{array}{cc}
1 & a \\ 
-a & 1%
\end{array}%
\right) ,$where $a=i\mathbf{i-j,}$ so that $a^{2}=0.$ This matrix is
self-dual with determinant $1$. It defines a random point field with exactly
2 points and correlation functions $R_{1}=1$ and $R_{2}=\left( 
\begin{array}{cc}
0 & 1 \\ 
1 & 0%
\end{array}%
\right) .$ On the other hand, $K$ does not have a diagonal form. (The
eigenvalues equal 1 but $K^{2}\neq K.$)

These examples shows that the conditions of Theorem \ref{theorem_existence}
are sufficient but not necessary. On the other hand, it is not possible to
drop entirely requirements on the eigenfunctions $u_{k}\left( x\right) $.
For example, let $\Lambda $ consist of two points and let $K=\frac{4}{3}%
\left( 
\begin{array}{cc}
1 & i/2 \\ 
i/2 & -1/4%
\end{array}%
\right) .$ (We consider this matrix as a quaternion matrix. It is self-dual
since $\left( i/2\right) ^{\ast }=i/2.$) This matrix has eigenvalues $1$ and 
$0$ and it can be computed that it has a diagonal form with $\lambda =1$ and 
$u=\left( 2/\sqrt{3}\right) [1,i/2]^{\ast }.$ However, it does not define a
valid random field since the first correlation function is negative at the
second point.

Many of the examples from random matrix theory are concerned with kernels
that do not have a real diagonal form. For example, one can check that the
circular orthogonal ensemble corresponds to a kernel without a diagonal form
and the circular symplectic ensemble corresponds to a kernel without a real
diagonal form. Hence, it appears desirable to find an extension of Theorem %
\ref{theorem_existence} . Ideally, we would like to know the sufficient and
necessary conditions that would ensure that a self-dual quaternion kernel
defines a valid random point field.

While this question is not answered in this paper, we can extend Theorem \ref%
{theorem_existence} and give sufficient conditions that ensure that a kernel
without a real diagonal form defines a valid random point field. First, let
us say that a kernel $K\left( x,y\right) $ has a \emph{quasi-real} \emph{%
diagonal form }(or simply call it \emph{quasi-real}) if it has a diagonal
form with real eigenvalues. Second, we will call it $\emph{positive}$ if 
\begin{equation*}
\mathrm{Det}_{M}\left. \left( K\left( x_{i},x_{j}\right) \right) \right\vert
_{1\leq i,j\leq m}\geq 0
\end{equation*}%
for all $m\in \mathbb{Z}^{+}$ and all $x_{1},\ldots ,x_{m}\in \Lambda $. (A
remarkable fact is that quasi-real kernels with positive eigenvalues are not
necessarily positive, as the last example above shows.)

We will call a quasi-real kernel $K\left( x,y\right) $ \emph{completely
positive} if every of the kernels 
\begin{equation*}
K_{I}=\sum_{i\in I}u_{i}\left( x\right) u_{i}^{\ast }\left( y\right)
\end{equation*}%
is positive, where $\{ u_{i}\left( x\right) \} $ is an orthonormal set of
eigenfunctions of $K\left( x,y\right) $ and $I=\left( i_{1},\ldots
,i_{m}\right) $ denote an \ ordered subset of indices of all eigenfunctions.

\begin{theorem}
\label{theorem_existence_quasireal} Suppose that a quaternion kernel $%
K\left( x,y\right) $ is finite-rank, quasi-real and completely positive.
Assume that all $\lambda _{k}\in \left[ 0,1\right] .$ Then the kernel $%
K\left( x,y\right) $ defines a Pfaffian point field.
\end{theorem}

Remark: We omit the assumption that $\sum_{k=1}^{\infty }\lambda _{k}<\infty 
$ which we imposed in Theorem \ref{theorem_existence} since we assume that
the kernel is finite-rank$.$ It should be possible to extend this theorem to
a more general case when there are infinite number of $\lambda _{k}$ and $%
\sum_{k=1}^{\infty }\lambda _{k}<\infty .$

Still, even if we restrict attention to kernels with a diagonal form, the
conditions of Theorem \ref{theorem_existence_quasireal} are not necessary.
For example, consider the two-point space $\Lambda $ with the counting
measure $\mu .$ Let $a=\left( 1+2i\right) +\left( \frac{19}{10}-\frac{20}{19}%
i\right) \mathbf{i,}$ and define the kernel as $K=\frac{1}{2}\left( 
\begin{array}{cc}
1 & a \\ 
a^{\ast } & 1%
\end{array}%
\right) .$ Then $\mathrm{Det}_{M}\left( K\right) \approx 0.3745,$ the pair
correlation function is positive, and the kernel defines a valid random
point field on $\Lambda $. On the other hand the eigenvalues are $\lambda
_{1,2}\approx \frac{1}{2}\pm 0.3529i$ and therefore the kernel is not
quasi-real.

In order to prove Theorems \ref{theorem_existence} and \ref%
{theorem_existence_quasireal}, let us introduce the following notations.

Suppose $K\left( x,y\right) $ has a diagonal form with eigenvalues $\lambda
_{k},$ that all $\lambda _{k}\in \left[ 0,1\right] ,$ and that $%
\sum_{k=1}^{\infty }\lambda _{k}<\infty .$ Define a random kernel 
\begin{equation}
K_{\xi }\left( x,y\right) =\sum_{k=1}^{\infty }\xi _{k}u_{k}\left( x\right)
u_{k}^{\ast }\left( y\right) ,  \label{def_q_kernel_random}
\end{equation}%
where $\xi _{k}$ are independent Bernoulli random variables. The random
variable $\xi _{k}$ takes value $1$ with probability $\lambda _{k}$. (We
will prove later that this kernel defines a valid point field.)

Theorems \ref{theorem_existence} and \ref{theorem_existence_quasireal} are
immediate consequences of the following results.

\begin{theorem}
\label{theorem_equivalence_Bernoulli}Suppose $K\left( x,y\right) $ has a
real diagonal form with eigenvalues $\lambda _{k},$ that all $\lambda
_{k}\in \left[ 0,1\right] ,$ and that $\sum_{k=1}^{\infty }\lambda
_{k}<\infty .$ Let $\mathcal{X}_{\xi }$ be a random point field which is a
mix of Pfaffian point fields with random kernels $K_{\xi }\left( x,y\right) $
defined in (\ref{def_q_kernel_random}). Then, $\mathcal{X}_{\xi }$ is a
Pfaffian point field with the kernel $K\left( x,y\right) $.
\end{theorem}

An analogous result holds for the quasi-real case.

\begin{theorem}
\label{theorem_equivalence_Bernoulli_quasi}Suppose $K\left( x,y\right) $ is
finite-rank, quasi-real and completely positive with eigenvalues $\lambda
_{k}\in \left[ 0,1\right] .$ Let $\mathcal{X}_{\xi }$ be a random point
field which is a mix of Pfaffian point fields with random kernels $K_{\xi
}\left( x,y\right) $ defined in (\ref{def_q_kernel_random}). Then, $\mathcal{%
X}_{\xi }$ is a Pfaffian point field with the kernel $K\left( x,y\right) $.
\end{theorem}

These theorems are quaternion analogues of Theorem 7 in \cite{hkpv06}. The
key ingredient in their proofs is an analogue of the Cauchy-Binet formula
for quaternion determinants which we state in Section \ref%
{section_Cauchy_Binet} and prove in Appendix. We will prove these theorems
in Section \ref{section_QDF}.

New kernels of Pfaffian point fields can also be obtained by the restriction
operation.

\begin{proposition}
Suppose $K\left( x,y\right) $ is a kernel of a $\Lambda $-valued Pfaffian
point field $\mathcal{X}$. Let $D\subset \Lambda .$ Then $K_{D}\left(
x,y\right) =\mathbf{1}_{D}\left( x\right) K\left( x,y\right) \mathbf{1}%
_{D}\left( y\right) $ is a kernel of another $\Lambda $-valued Pfaffian
point field.
\end{proposition}

\textbf{Proof:} The correlation functions defined by the kernel $K_{D}\left(
x,y\right) $ are valid correlation functions since they equal the
correlations functions of a random point field $\mathcal{X}_{D},$ generated
by the following procedure. First, generate the points of $\mathcal{X}$.
Then remove the points which are outside of $D.$ $\square $

In addition to the existence results, Theorem \ref%
{theorem_equivalence_Bernoulli} allows us to study the total number of
points in a Pfaffian point field.

\begin{theorem}
\label{theorem_characteristic_function}Let $\mathcal{X}$ be a $\Lambda $%
-valued Pfaffian point field with a finite-rank kernel $K.$ Let $\mathcal{N}$
denote the number of points of $\mathcal{X}$ in $\Lambda .$ Then the
characteristic function of the random variable $\mathcal{N}$ satisfies the
following equations:%
\begin{eqnarray}
\varphi _{\mathcal{N}}\left( t\right) &\equiv &\mathbb{E}\left( e^{i\mathcal{%
N}t}\right) =\prod\limits_{k=1}^{r}\left( 1+\left( e^{it}-1\right) \lambda
_{k}\right)  \notag \\
&=&\mathrm{Det}_{M}\left( I+\left( e^{it}-1\right) K\right) ,
\label{char_function_of_N}
\end{eqnarray}%
where $\lambda _{k}$ are eigenvalues of the kernel $K$ and $r$ is its rank.
\end{theorem}

A consequence of this result is the central limit theorem for the number of
points in Pfaffian point fields.

\begin{theorem}
\label{theorem_CLT_general}Let $\mathcal{X}_{n}$ be a sequence of $\Lambda
_{n}$-valued Pfaffian point fields with finite-rank kernels $K_{n}$. Let $%
\mathcal{N}_{n}$ denote the number of points of $\mathcal{X}_{n}.$ Suppose
that all eigenvalues of kernels $K_{n}$ are real and in the interval $\left[
0,1\right] ,$ and that $\mathbb{V}ar\left( \mathcal{N}_{n}\right)
\rightarrow \infty $ as $n\rightarrow \infty .$ Then, the sequence of random
variables $\left( \mathcal{N}_{n}-\mathbb{E}\left( \mathcal{N}_{n}\right)
\right) /\sqrt{\mathbb{V}ar\left( \mathcal{N}_{n}\right) }$ approaches a
standard Gaussian random variable in distribution.
\end{theorem}

This is an analog of a theorem that was proved by Costin and Lebowitz in 
\cite{costin_lebowitz95} and Diaconis and Evans in \cite{diaconis_evans01}
for particular cases, and by Soshnikov in \cite{soshnikov00a} for general
determinantal ensembles. Later, a simplified proof was suggested in \cite%
{hkpv06} and we use its main idea to prove our theorem.

We will see in the final section that the circular and Gaussian symplectic
ensemble of random matrices have finite-rank projection kernels . Numerical
evaluations suggest that the restrictions of these kernels have real
eigenvalues in the interval $\left[ 0,1\right] $. However, the proof of this
claim is elusive.

In the theory of determinantal random fields, a prominent place is given to
determinantal formulas for evaluation of expressions like $\mathbb{E}%
\prod\limits_{k=1}^{N}\left( 1+f\left( x_{k}\right) \right) ,$ where $x_{k}$
denote the points of the field. These formulas are useful for calculating
the distribution of spacings and similar quantities. It turns out that
similar expressions can be written for Pfaffian fields. (This was also
observed in Rains \cite{rains00} in a somewhat different form.)

First, we define the quaternion version of the Fredholm determinant for the
self-dual kernel $K\left( x,y\right) $:%
\begin{equation*}
\mathrm{Det}_{M}\left( I+K\right) :=1+\sum_{n=1}^{\infty }\frac{1}{n!}%
\int_{\Lambda ^{n}}\mathrm{Det}_{M}\left. \left( K\left( x_{i},x_{j}\right)
\right) \right\vert _{1\leq i,j\leq n}dx_{1}\ldots dx_{n}.
\end{equation*}%
Note that for finite-rank kernels the series has only finite number of
terms. In addition, it can be checked that if $\Lambda $ is finite, so that $%
I+K$ is a matrix, then this definition agrees with the usual definition of
the Dyson-Moore determinant for the matrix $I+K.$

Let us consider $\mathbb{E}\prod_{k}\left( 1+f\left( x_{k}\right) \right) ,$
where $f\left( x\right) $ is a complex-valued function and the product
extended over all points of the field.

\begin{theorem}
\label{prop_Fredholm_determinant_formula}Suppose that the kernel $K$ of a
Pfaffian field has finite rank. Then%
\begin{equation*}
\mathbb{E}\prod_{k}\left( 1+f\left( x_{k}\right) \right) =\mathrm{Det}%
_{M}\left( I+\sqrt{f\left( x\right) }K\sqrt{f\left( y\right) }\right) ,
\end{equation*}%
whenever terms on both sides of this equality are well-defined.
\end{theorem}

\textbf{Proof:} The assumption implies that the field has a finite number of
points and the product $\mathbb{E}\prod_{k}\left( 1+f\left( x_{k}\right)
\right) $ has finite number of terms. By the definition of correlation
functions%
\begin{eqnarray*}
\mathbb{E}\prod_{i}\left( 1+f\left( x_{i}\right) \right)
&=&1+\sum_{n=1}^{\infty }\frac{1}{n!}\int_{\Lambda ^{n}}\left(
\prod_{i=1}^{n}f\left( x_{i}\right) \right) R\left( x_{1},\ldots
,x_{n}\right) dx_{1}\ldots dx_{n} \\
&=&1+\sum_{n=1}^{\infty }\frac{1}{n!}\int_{\Lambda ^{n}}\left(
\prod_{i=1}^{n}f\left( x_{i}\right) \right) \mathrm{Det}_{M}\left. \left(
K\left( x_{i},x_{j}\right) \right) \right\vert _{1\leq i,j\leq
n}dx_{1}\ldots dx_{n} \\
&=&1+\sum_{n=1}^{\infty }\frac{1}{n!}\int_{\Lambda ^{n}}\mathrm{Det}%
_{M}\left. \left( \sqrt{f\left( x_{i}\right) }K\left( x_{i},x_{j}\right) 
\sqrt{f\left( x_{j}\right) }\right) \right\vert _{1\leq i,j\leq
n}dx_{1}\ldots dx_{n} \\
&=&\mathrm{Det}_{M}\left( I+\sqrt{f\left( x\right) }K\sqrt{f\left( y\right) }%
\right) .
\end{eqnarray*}%
(Note that in fact all the series in this calculation have the finite number
of terms.) $\square $

In particular, formula (\ref{char_function_of_N}) can be alternatively
obtained from Theorem \ref{prop_Fredholm_determinant_formula} by taking $%
f\left( x\right) =e^{it}-1,$ a function which is constant in $x$ and depends
only on a parameter $t.$

The rest of the paper is organized as follows. Section \ref%
{section_Cauchy_Binet} formulates a quaternion version of the Cauchy-Binet
identity. Section \ref{section_QDF} proves the existence Theorems \ref%
{theorem_existence} and \ref{theorem_existence_quasireal}, Section \ref%
{section_CLT} proves the CLT in Theorem \ref{theorem_CLT_general}, and
Section \ref{section_example} provides an illustration by considering the
circular and Gaussian symplectic ensembles of random matrices. The appendix
contains background information about quaternion matrices and determinants
and a proof of the quaternion Cauchy-Binet identity.

\section{Cauchy-Binet formula}

\label{section_Cauchy_Binet}

Determinants of quaternion matrices have been studied for almost a hundred
years. (Study wrote a paper \cite{study20} \ about them in 1920, and Moore
presented his definition \cite{moore22} in 1922). In the 1970s, this subject
got a boost after Dyson (\cite{dyson70}) re-discovered Moore's determinant
and related it to the distribution of eigenvalues of random matrices.

Unfortunately, none of the available quaternion determinants enjoys all the
properties of the usual determinant, and the validity of each standard
determinantal identity has to be checked individually.

The Cauchy-Binet identity states \ that if $A$ is an $m$-by-$n$ matrix and $%
B $ is an $n$-by-$m$ matrix, with $n\geq m,$ then 
\begin{equation*}
\det \left( AB\right) =\sum_{I}\det \left( A^{I}\right) \det \left(
B^{I}\right) ,
\end{equation*}%
where the summation is over $I=\left( i_{1}<i_{2}<\ldots <i_{m}\right) $,
the ordered subsets of $\left\{ 1,\ldots ,n\right\} $ that consist of $m$
elements. Matrices $A^{I}$ and $B^{I}$ are square matrices that consist of $%
m $ columns of $A$ and $m$ rows of $B,$ respectively, with indices in $I.$

An implicit assumption in this result is that the entries of the matrices
\thinspace $A$ and $B$ are from a commutative ring, for example from a field
of complex numbers. Unfortunately, in this form the Cauchy-Binet identity
fails for the quaternion determinants.

However, a weaker form of the Cauchy-Binet identity still holds.

\begin{theorem}
\label{theorem_cauchy_binet} Let $C$ be an $n$-by-$m$ quaternion matrix, $%
n\geq m,$ and let $C^{\ast }$ be the dual of $C$. Then 
\begin{equation*}
\mathrm{Det}_{M}\left( C^{\ast }C\right) =\sum_{I}\mathrm{Det}_{M}\left(
\left( C^{I}\right) ^{\ast }C^{I}\right) ,
\end{equation*}%
where the summation is over $I=\left( i_{1}<i_{2}<\ldots <i_{m}\right) $,
the ordered subsets of $\left\{ 1,\ldots ,n\right\} $ that consist of $m$
elements, and where $C^{I}$ consists of rows $i_{1}$,\ldots ,$i_{m}$ of $C.$
\end{theorem}

A proof of this theorem is in Appendix.

\begin{corollary}
\label{corollary_Cauchy_Binet}Suppose that $\Lambda $ is an $n$-by-$n$
diagonal matrix with scalar entries $\lambda _{i}$ on the main diagonal and
that $C$ is an $n$-by-$m$ quaternion matrix, $n\geq m$. Then, we have. 
\begin{equation*}
\mathrm{Det}_{M}\left( C^{\ast }\Lambda C\right) =\sum_{\substack{ I=\left(
i_{1},\ldots ,i_{m}\right)  \\ i_{1}<\ldots <i_{m}}}\lambda _{i_{1}}\ldots
\lambda _{i_{m}}\mathrm{Det}_{M}\left( \left( C^{I}\right) ^{\ast
}C^{I}\right) .
\end{equation*}
\end{corollary}

\textbf{Proof of Corollary \ref{corollary_Cauchy_Binet}}: \ Let $\Lambda
^{1/2}$ be an $n$-by-$n$ diagonal matrix with scalar entries such that $%
\left( \Lambda ^{1/2}\right) ^{2}=\Lambda .$ For $I=\left( i_{1},\ldots
,i_{m}\right) ,$ let $\left( \Lambda ^{1/2}\right) ^{II}$ denote an $m$-by-$%
m $ matrix which is formed by taking entries at the intersection of rows and
columns $i_{1},$ \ldots , $i_{m}$ in matrix $\Lambda ^{1/2}.$ We write 
\begin{eqnarray*}
\mathrm{Det}_{M}\left( C^{\ast }\Lambda C\right) &=&\sum_{I}\mathrm{Det}%
_{M}\left( \left( \Lambda ^{1/2}C\right) ^{I\ast }\left( \Lambda
^{1/2}C\right) ^{I}\right) \\
&=&\sum_{I}\mathrm{Det}_{S}\left( \left( \Lambda ^{1/2}C\right) ^{I}\right)
\\
&=&\sum_{I}\mathrm{Det}_{S}\left( \left( \Lambda ^{1/2}\right)
^{II}C^{I}\right) \\
&=&\sum_{I}\mathrm{Det}_{S}\left( \left( \Lambda ^{1/2}\right) ^{II}\right) 
\mathrm{Det}_{S}\left( C^{I}\right) \\
&=&\sum_{I=\left( i_{1},\ldots ,i_{m}\right) }\lambda _{i_{1}}\ldots \lambda
_{i_{m}}\mathrm{Det}_{M}\left( \left( C^{I}\right) ^{\ast }C^{I}\right) .
\end{eqnarray*}

The first line is the Cauchy-Binet identity. The second line is the relation
between the Moore-Dyson and Study determinants (see (\ref%
{identity_study_moore}) in Appendix). The fourth line follows by
multiplicativity of the Study determinant. And in the fifth line we have
used (\ref{identity_study_moore}) again. $\square $

\section{The existence of Pfaffian point fields}

\label{section_QDF}

\begin{proposition}
\label{prop_projection_kernel}Let $K_{N}\left( x,y\right)
=\sum_{k=1}^{N}u_{k}\left( x\right) u_{k}^{\ast }\left( y\right) ,$ where $%
u_{k}\left( x\right) $ are orthonormal quaternion functions, and assume that 
$K_{N}\left( x,y\right) $ is positive, that is, 
\begin{equation*}
\mathrm{Det}_{M}\left. \left( K_{N}\left( x_{i},x_{j}\right) \right)
\right\vert _{1\leq i,j\leq m}\geq 0
\end{equation*}%
for all $m\in \mathbb{Z}^{+}$ and all $x_{1},\ldots ,x_{m}$. Then $%
K_{N}\left( x,y\right) $ defines a valid symplectic determinantal field with
exactly $N$ points.
\end{proposition}

Note that for the case when $u_{k}\left( x\right) $ are real quaternion
functions the assumption of positivity can be dropped. Indeed, in this case,
the matrix $K=\left. K_{r}\left( x_{i},x_{j}\right) \right\vert
_{i,j=1,\ldots ,m}$ is positive semidefinite. (That is, for every real
quaternion vector $v,$ $v^{\ast }Kv\geq 0.$) This implies that all
eigenvalues of this matrix are real and non-negative. For self-dual matrices
with real quaternion entries the Moore-Dyson determinant can be computed as
the product of the eigenvalues, and we conclude that $\mathrm{Det}_{M}\left.
\left( K_{N}\left( x_{i},x_{j}\right) \right) \right\vert _{1\leq i,j\leq
m}\geq 0$.

\textbf{Proof of Proposition \ref{prop_projection_kernel}:} We need to show
that the functions defined by the rule 
\begin{equation*}
R_{m}\left( x_{1},\ldots x_{m}\right) =\mathrm{Det}_{M}\left( K_{N}\left(
x_{i},x_{j}\right) \right) ,\text{ }i,j=1,\ldots ,m
\end{equation*}%
are the correlation functions of a random field.

First, by assumption of positivity all functions $R_{m}\left( x_{1},\ldots
x_{m}\right) $ are non-negative.

Next, by integrating the kernel we find that 
\begin{eqnarray*}
\int_{\mathbb{R}}K_{N}\left( x,x\right) dx &=&N, \\
\int_{\mathbb{R}}K_{N}\left( x,y\right) K_{N}\left( y,x\right) dydx &=&N.
\end{eqnarray*}%
By using formulas (\ref{def_correlation_formulas}), we can conclude that the
total number of points in the process with kernel $K_{N}\left( x,y\right) $
is exactly $N$. (Its expectation is $N$ and its variance is $0.$) Hence, it
remains to show that the functions $R_{m}\left( x_{1},\ldots x_{m}\right) $
agree among themselves for all $m$. This can be done by the quaternion
analogue of the Mehta lemma. (Compare Theorem 5.1.4 on p. 75 in Mehta \cite%
{mehta04})

Let $K_{m}:=(K(x_{i},x_{j}))_{1\leq i,j\leq m}$ for a self-dual quaternion
kernel $K\left( x,y\right) $.

\begin{lemma}[Dyson]
Assume that $K$ satisfies either 
\begin{equation}
\int_{\mathbb{R}}K(x,y)K(y,z)dy=K(x,z)  \label{kernel_formula2}
\end{equation}%
or 
\begin{equation}
\varphi \left( \int_{\mathbb{R}}K(x,y)K(y,z)dy\right) =\varphi \left(
K(x,z)\right) +E\varphi \left( K(x,z)\right) -\varphi \left( K(x,z)\right) E,
\label{kernel_formula3}
\end{equation}%
where 
\begin{equation}
E=\left( 
\begin{array}{cc}
1 & 0 \\ 
0 & 0%
\end{array}%
\right) .
\end{equation}%
Then 
\begin{equation}
\int_{\mathbb{R}}\mathrm{Det}_{M}(K_{m})dx_{m}=(N-m+1)\mathrm{Det}%
_{M}(K_{m-1}),  \label{Mehta_formula}
\end{equation}%
where $N=\int K(x,x)d\mu (x)$.
\end{lemma}

This result is due to Dyson. (See proof of Theorem 4 in \cite{dyson70}.)

For the kernel $K_{N}\left( x,y\right) ,$ equation (\ref{kernel_formula2})
holds, and we find that 
\begin{equation*}
\int_{\mathbb{R}}R_{m}\left( x_{1},\ldots ,x_{m}\right) dx_{m}=\left(
N-(m-1)\right) R_{m-1}\left( x_{1},\ldots ,x_{m-1}\right) ,
\end{equation*}%
which shows that all correlation functions are all in agreement.

\begin{corollary}
\label{prop_existence_Bernoulli}Assume that 
\begin{equation}
K_{\xi }\left( x,y\right) =\sum_{k=1}^{\infty }\xi _{k}u_{k}\left( x\right)
u_{k}^{\ast }\left( y\right) ,  \label{def_q_kernel_random_2}
\end{equation}%
where every $\xi _{k}$ is an independent Bernoulli random variable that
takes value $1$ with probability $\lambda _{k}$, and where $u_{k}\left(
x\right) $ are orthonormal real quaternion functions. Suppose that $%
\sum_{k=1}^{\infty }\lambda _{k}<\infty $. Conditional on $\xi $, the
function $K_{\xi }\left( x,y\right) $ is a kernel of a determinantal point
field, $\mathcal{X}_{\xi },$ and the total number of points in $\mathcal{X}%
_{\xi }$ is $\sum_{k=1}^{\infty }\xi _{k}.$
\end{corollary}

We can prove an analogous result for a quasi-real diagonal form.

\begin{corollary}
\label{prop_existence_Bernoulli_quasi}Assume that 
\begin{equation}
K_{\xi }\left( x,y\right) =\sum_{k=1}^{N}\xi _{k}u_{k}\left( x\right)
u_{k}^{\ast }\left( y\right) ,
\end{equation}%
where every $\xi _{k}$ is an independent Bernoulli random variable that
takes value $1$ with probability $\lambda _{k}$. Suppose that the kernel $%
K\left( x,y\right) =\sum_{k=1}^{N}u_{k}\left( x\right) u_{k}^{\ast }\left(
y\right) $ is completely positive. Conditional on $\xi $, the function $%
K_{\xi }\left( x,y\right) $ is a kernel of a determinantal point field, $%
\mathcal{X}_{\xi }$ and the total number of points in $\mathcal{X}_{\xi }$
is $\sum_{k=1}^{N}\xi _{k}.$
\end{corollary}

\textbf{Proof of Theorem \ref{theorem_equivalence_Bernoulli}:} By the
iterated expectation formula, the correlation functions of process $\mathcal{%
X}_{\xi }$ equal $\mathbb{E}\mathrm{Det}_{M}\left( K_{\xi }\left(
x_{i},x_{j}\right) \right) ,$ where the expectation is taken over randomness
in $\xi $. Hence, it is enough to prove that

\begin{equation}
\mathbb{E}\mathrm{Det}_{M}\left( K_{\xi }\left( x_{i},x_{j}\right) \right) =%
\mathrm{Det}_{M}\left( K\left( x_{i},x_{j}\right) \right) ,\text{ }%
i,j=1,\ldots ,m,\text{ }  \label{determinant_equality}
\end{equation}%
almost everywhere.

First, let 
\begin{equation*}
K_{R}\left( x,y\right) =\sum_{k=1}^{R}\lambda _{k}u_{k}\left( x\right)
u_{k}^{\ast }\left( y\right) .
\end{equation*}%
Since we assumed the absolute convergence of the kernel, hence $\mathrm{Det}%
_{M}\left( K_{R}\left( x_{i},x_{j}\right) \right) $ converges to $\mathrm{Det%
}_{M}\left( K\left( x_{i},x_{j}\right) \right) $ almost everywhere as $%
R\rightarrow \infty .$

Let $R$-by-$m$ matrix $C$ be defined as 
\begin{equation*}
C_{kl}=u_{k}^{\ast }\left( x_{l}\right) ,\text{ }k=1,\ldots ,R;\text{ }%
l=1,\ldots ,m.
\end{equation*}%
and let $\Lambda $ be an $R$-by-$R$ diagonal matrix with diagonal entries $%
\lambda _{i}$.

By Corollary \ref{corollary_Cauchy_Binet}, 
\begin{equation}
\mathrm{Det}_{M}\left( K_{R}\left( x_{i},x_{j}\right) \right) =\sum 
_{\substack{ I=\left( i_{1},\ldots ,i_{m}\right)  \\ i_{1}<\ldots <i_{m}}}%
\lambda _{i_{1}}\ldots \lambda _{i_{m}}\mathrm{Det}_{M}\left( \left(
C^{I}\right) ^{\ast }C^{I}\right) .  \label{determinant_original_field}
\end{equation}

Next, let the random variable $\mathrm{Det}_{M}\left( K_{\xi }\left(
x_{i},x_{j}\right) \right) $ be denoted as $Y$ and let $A_{R}$ be the event
that all $\xi _{k}$ are zero for $k>R.$ (That is, $A_{R}=\cap _{k>R}\{\xi
_{k}=0\}$). Note that%
\begin{equation}
\mathbb{E}Y=\mathbb{E}\left( Y|A_{R}\right) \mathbb{P}\left( A_{R}\}\right) +%
\mathbb{E}\left( Y|A_{R}^{c}\right) \mathbb{P}\left( A_{R}^{c}\right) .
\label{equality_conditional_on_R}
\end{equation}%
By using independence of $A_{R}$ and $\xi _{k}$ for $k\leq R,$ we find that 
\begin{equation*}
\mathbb{E}\left( Y|A_{R}\right) \mathbb{P}\left( A_{R}\right) =\mathbb{E}%
\mathrm{Det}_{M}\left( C^{\ast }\Lambda _{\xi }C\right) \mathbb{P}\left(
A_{R}\right) ,
\end{equation*}%
where $\Lambda _{\xi }$ denotes an $R$-by-$R$ diagonal matrix with diagonal
entries $\xi _{i}$.

Next, by Corollary \ref{corollary_Cauchy_Binet}, 
\begin{equation*}
\mathbb{E}\mathrm{Det}_{M}\left( C^{\ast }\Lambda _{\xi }C\right) =\mathbb{E}%
\sum_{\substack{ I=\left( i_{1},\ldots ,i_{m}\right)  \\ i_{1}<\ldots <i_{m} 
}}\xi _{i_{1}}\ldots \xi _{i_{m}}\mathrm{Det}_{M}\left( \left( C^{I}\right)
^{\ast }C^{I}\right) .
\end{equation*}

Since the variables $\xi _{i_{1}},\ldots ,\xi _{i_{m}}$ are independent and
have expectation $\lambda _{i_{1}},\ldots ,\lambda _{i_{m}},$ we find that $%
\mathbb{E}\xi _{i_{1}}\ldots \xi _{i_{m}}=\lambda _{i_{1}}\ldots \lambda
_{i_{m}}.$ Hence 
\begin{equation}
\mathbb{E}\left( Y|A_{R}\}\right) \mathbb{P}\left( A_{R}\right) =\mathrm{Det}%
_{M}\left( K_{R}\left( x_{i},x_{j}\right) \right) \mathbb{P}\left(
A_{R}\right) ,  \label{determinant_random_field}
\end{equation}%
and the probability $\mathbb{P}\left( A_{R}\right) $ converges to $1$ as $%
R\rightarrow \infty $ by independence of $\xi _{k},$ Borel-Cantelli lemma
and the assumption $\sum \lambda _{k}<\infty .$

Now let us show that $\mathbb{E}\left( Y|A_{R}^{c}\}\right) \mathbb{P}\left(
A_{R}^{c}\right) $ converges to zero almost everywhere. By positivity of the
determinant, it is enough to show that 
\begin{equation}
\int_{\mathbb{R}^{m}}\mathbb{E}\left( \mathrm{Det}_{M}\left( K_{\xi }\left(
x_{i},x_{j}\right) \right) |A_{R}^{c}\right) \mathbb{P}\left(
A_{R}^{c}\right) \rightarrow 0  \label{condition_convergence}
\end{equation}%
as $R\rightarrow \infty .$

Let 
\begin{equation*}
n_{\xi }:=\sum_{k=1}^{\infty }\xi _{k},
\end{equation*}%
which is finite since both the expectation and the variance of the sum on
the right hand side are convergent. Since $K_{\xi }$ is a projection
operator, the total number of points of the process $\mathcal{X}_{\xi }$ in $%
\mathbb{R}$ equals $n_{\xi }.$ By changing the order of the expectation and
the integral signs in (\ref{condition_convergence}), which is possible since
the integrand is positive, and by using the identities for correlation
functions we obtain that we need to estimate%
\begin{equation*}
\mathbb{E}\left( n_{\xi }\left( n_{\xi }-1\right) \ldots \left( n_{\xi
}-m+1\right) |A_{R}^{c}\right) \mathbb{P}\left( A_{R}^{c}\right) ,
\end{equation*}%
which is smaller than 
\begin{equation*}
\mathbb{E}\left( n_{\xi }^{m}|A_{R}^{c}\right) \mathbb{P}\left(
A_{R}^{c}\right) .
\end{equation*}%
By expanding 
\begin{equation*}
n_{\xi }^{m}=\left( \sum_{k=1}^{\infty }\xi _{k}\right) ^{m},
\end{equation*}%
and using the fact that $\xi _{k}^{s}=\xi _{k}$ for every integer $s\geq 1,$
we observe that it is enough to show that 
\begin{equation*}
\mathbb{E}\left( \sum_{i_{1}<\ldots <i_{r}}\xi _{i_{1}}\xi _{i_{2}}\ldots
\xi _{i_{r}}|A_{R}^{c}\right) \mathbb{P}\left( A_{R}^{c}\right) \rightarrow
0,
\end{equation*}%
as $R\rightarrow \infty ,$ where the sum is over all possible ordered $r$%
-tuples $\left( i_{1},\ldots ,i_{r}\right) $ such that $i_{1}<\ldots <i_{r}$%
, and $1\leq r\leq m.$

We can divide the sum in two parts. The first part is when $i_{r}\leq R.$ In
this case, the variables $\xi _{i_{1}}$, $\xi _{i_{2}}$, $\ldots ,$ $\xi
_{i_{r}}$ are independent from the event $A_{R}^{c},$ and therefore we can
estimate this part of the sum as 
\begin{equation*}
\left( \sum_{i_{1}<\ldots <i_{r}}\lambda _{i_{1}}\lambda _{i_{2}}\ldots
\lambda _{i_{r}}|A_{R}^{c}\right) \mathbb{P}\left( A_{R}^{c}\right) \leq
cS^{m}\mathbb{P}\left( A_{R}^{c}\right) ,
\end{equation*}%
where $S:=\sum_{k=1}^{\infty }\lambda _{k}$ and $c$ is an absolute constant.
Hence this part converges to zero as $R\rightarrow \infty $ because $\mathbb{%
P}\left( A_{R}^{c}\right) $ converges to zero.

The second part of the sum is when $i_{r}>R.$ In this case, the event $\xi
_{i_{r}}=1$ implies $A_{R}^{c}$ and therefore, 
\begin{eqnarray*}
\mathbb{E}\left( \xi _{i_{1}}\xi _{i_{2}}\ldots \xi
_{i_{r}}|A_{R}^{c}\right) \mathbb{P}\left( A_{R}^{c}\right) &=&\mathbb{P}%
\left( \xi _{i_{1}}=1,\ldots ,\xi _{i_{r}}=1,\text{ and }A_{R}^{c}\right) \\
&=&\mathbb{P}\left( \xi _{i_{1}}=1,\ldots ,\xi _{i_{r}}=1\right) \\
&=&\lambda _{i_{1}}\ldots \lambda _{i_{r}}.
\end{eqnarray*}%
Therefore, we estimate the second part as 
\begin{eqnarray*}
\sum_{\substack{ i_{1}<\ldots <i_{r},  \\ i_{r}>R}}\lambda _{i_{1}}\ldots
\lambda _{i_{r}} &\leq &\sum_{i_{r}=R+1}^{\infty }\sum_{i_{1},\ldots
i_{r-1}=1}^{\infty }\lambda _{i_{1}}\ldots \lambda _{i_{r}} \\
&=&S^{m-1}\sum_{i_{r}=R+1}^{\infty }\lambda _{i_{r}}.
\end{eqnarray*}

Hence the second part converges to zero as $R\rightarrow \infty ,$ because
the tail sums $\sum_{i_{r}=R+1}^{\infty }\lambda _{i_{r}}$ converge to zero.

This shows that $\mathbb{E}\left( Y|A_{R}^{c}\}\right) \mathbb{P}\left(
A_{R}^{c}\right) $ in (\ref{equality_conditional_on_R}) converges to zero as 
$R\rightarrow \infty .$ If we compare (\ref{determinant_original_field}) and
(\ref{determinant_random_field}) and let $R$ grow to infinity, we find that 
\begin{equation}
\mathbb{E}\mathrm{Det}_{M}\left( K_{\xi }\left( x_{i},x_{j}\right) \right) =%
\mathrm{Det}_{M}\left( K\left( x_{i},x_{j}\right) \right)
\end{equation}%
almost everywhere, and this completes the proof of the theorem. $\square .$

The proof of Theorem \ref{theorem_equivalence_Bernoulli_quasi}\textbf{\ }is
essentially the same as for Theorem \ref{theorem_equivalence_Bernoulli} but
it is easier, since we assumed that the number of eigenvalues is finite.

\section{Number of points in Pfaffian point fields}

\label{section_CLT}

\textbf{Proof of Theorem \ref{theorem_characteristic_function}:} The moments
of the number of points distribution can be written as polynomials in traces
of the kernel and its powers. Therefore, they can be written as symmetric
polynomials of the kernel eigenvalues. It follows that the characteristic
function of $\mathcal{N}$ can also be written as a symmetric function of the
kernel eigenvalues. The form of this function can be recovered from the
particular case when the kernel has a real diagonal form and Theorem \ref%
{theorem_equivalence_Bernoulli} is applicable. In this case $\mathcal{N}$ is
the sum of independent Bernoulli random variables, $\mathcal{N}=\xi
_{1}+\ldots +\xi _{r},$ where $\xi _{k}$ equals $0$ or $1$ with
probabilities $1-\lambda _{k}$ and $\lambda _{k},$ respectively. By
properties of characteristic functions, 
\begin{equation*}
\varphi _{\mathcal{N}}\left( t\right) =\prod\limits_{k=1}^{r}\varphi _{%
\mathcal{\xi }_{k}}\left( t\right) =\prod\limits_{k=1}^{r}\left( 1+\left(
e^{it}-1\right) \lambda _{k}\right) .
\end{equation*}%
The second equality follows because the determinant can be written as the
product of eigenvalues. $\square $

\textbf{Proof of Theorem \ref{theorem_CLT_general}: }The expression for the
characteristic function of $\mathcal{N}_{n}$ is 
\begin{equation*}
\varphi _{\mathcal{N}_{n}}\left( t\right) =\prod\limits_{k=1}^{r_{n}}\left(
1+\left( e^{it}-1\right) \lambda _{k}^{\left( n\right) }\right) ,
\end{equation*}%
where $\lambda _{k}^{\left( n\right) }$ are eigenvalues of the kernel $%
K_{n}. $ This expression is the same as for the sum of independent Bernoulli
random variables $\xi _{k}^{\left( n\right) }$ with $\mathbb{P}\left( \xi
_{k}^{\left( n\right) }=1\right) =\lambda _{k}^{\left( n\right) }.$ By the
Lindenberg-Feller theorem (\cite{shiryaev96}, Theorem III.4.2 on p. 334) we
conclude that if $\mathbb{V}ar\left( \mathcal{N}_{n}\right) \rightarrow
\infty $ as $n\rightarrow \infty ,$ then 
\begin{equation*}
\frac{\mathcal{N}_{n}-\mathbb{E}\left( \mathcal{N}_{n}\right) }{\sqrt{%
\mathbb{V}ar\left( \mathcal{N}_{n}\right) }}
\end{equation*}%
approaches the standard Gaussian random variable. $\square $

\section{Example: symplectic ensembles of random matrices}

\label{section_example}

\subsection{Circular Symplectic Ensemble}

The circular symplectic ensemble of random matrices (CUE) is defined as the
probability space of $N$-by-$N$ self-dual real quaternion unitary matrices.
The probability space has the Haar measure.

The eigenvalues of matrices from this ensembles are located on the unit
circle and can be identified with angles $\theta _{k}.$ The density of the
eigenvalue distribution is%
\begin{equation*}
c\prod_{1\leq j<k\leq N}\left\vert e^{i\theta _{j}}-e^{i\theta
_{k}}\right\vert ^{4},\text{ }-\pi \leq \theta _{l}<\pi .
\end{equation*}
(For more details see Dyson's papers or Chapter 2 in Forrester's book \cite%
{forrester10}.)

Let $N$ be any positive integer. Define 
\begin{equation*}
s_{2N}\left( \theta \right) :=\frac{1}{2\pi }\sum_{p}e^{ip\theta }=\frac{1}{%
\pi }\sum_{p>0}\cos \left( p\theta \right) =\frac{1}{2\pi }\frac{\sin \left(
N\theta /2\right) }{\sin \left( \theta /2\right) }.
\end{equation*}%
(Here the first summation is over $p=\left( -2N+1\right) /2,$ $\left(
-2N+3\right) /2,$ \ldots , $\left( 2N-1\right) /2,$ and the second summation
is over $p=1/2,\ldots ,\left( 2N-1\right) /2.$) Note that $s_{2N}\left(
\theta \right) $ is even in $\theta .$

Following Dyson, we write 
\begin{equation*}
Ds_{2N}\left( \theta \right) :=\frac{d}{d\theta }s_{2N}\left( \theta \right)
=\frac{i}{2\pi }\sum_{p}pe^{ip\theta }=-\frac{1}{\pi }\sum_{p>0}p\sin \left(
p\theta \right) ,
\end{equation*}%
and 
\begin{equation*}
Is_{2N}\left( \theta \right) :=\int_{0}^{\theta }s_{2N}\left( \theta
^{\prime }\right) d\theta ^{\prime },
\end{equation*}%
so that 
\begin{equation*}
Is_{2N}\left( \theta \right) =\frac{1}{2\pi i}\sum_{p}p^{-1}e^{ip\theta }=%
\frac{1}{\pi }\sum_{p>0}\frac{1}{p}\sin \left( p\theta \right)
\end{equation*}%
The functions $Ds_{2N},$ and $Is_{2N}$ are odd in $\theta .$

Define the quaternion function $\sigma _{N4}\left( \theta \right) $ by its
matrix representation: 
\begin{equation*}
\varphi \left( \sigma _{N4}\left( \theta \right) \right) =\frac{1}{2}\left( 
\begin{array}{cc}
s_{2N}\left( \theta \right) & Ds_{2N}\left( \theta \right) \\ 
Is_{2N}\left( \theta \right) & s_{2N}\left( \theta \right)%
\end{array}%
\right) =\frac{1}{2\pi }\sum_{p>0}\left( 
\begin{array}{cc}
\cos \left( p\theta \right) & -p\sin \left( p\theta \right) \\ 
p^{-1}\sin \left( p\theta \right) & \cos \left( p\theta \right)%
\end{array}%
\right) .
\end{equation*}%
In terms of quaternions, the kernel can be written as follows:

\begin{eqnarray*}
\sigma _{N4}\left( \theta \right) &=&\frac{1}{2}\left( s_{2N}-\frac{1}{2}%
\left( Is_{2N}+Ds_{2N}\right) i\mathbf{i+}\frac{1}{2}\left(
Is_{2N}-Ds_{2N}\right) \mathbf{j}\right) \\
&=&\frac{1}{2\pi }\sum_{p=1/2}^{N-1/2}(\cos p\theta +a_{p}\sin p\theta ),
\end{eqnarray*}%
where 
\begin{equation*}
a_{p}=\frac{1}{2p}\left[ \left( p^{2}-1\right) i\mathbf{i+}\left(
p^{2}+1\right) \mathbf{j}\right] .
\end{equation*}%
It is easy to check that $a_{p}^{2}=-1.$

Dyson proved the following result (See Theorem 3 in \cite{dyson70}): The
random field of eigenvalues is Pfaffian with the kernel $\sigma _{N4}\left(
\theta -\theta ^{\prime }\right) $.

The kernel $\sigma _{N4}$ is a projection on a subspace of a
finite-dimensional linear space (or rather module) $L$ over $\mathbb{Q}.$ 
\begin{equation*}
L=\mathrm{span}\left\{ \frac{1}{\sqrt{\pi }}\cos p\theta ,\text{ }\frac{1}{%
\sqrt{\pi }}\sin p\theta ,\text{ }p=\frac{1}{2},\frac{3}{2},\ldots ,\frac{%
2N-1}{2}\right\} .
\end{equation*}%
Note that $\dim _{\mathbb{Q}}L=2N.$

If we calculate the action of the kernel $\sigma _{N4}\left( \theta \right) $
in the basis $\pi ^{-1/2}\cos p\theta ,$ $\pi ^{-1/2}\sin p\theta ,$ then we
find the following matrix representation of the operator with kernel $\sigma
_{N4}:$ 
\begin{equation}
K=\frac{1}{2}\left\{ \left( 
\begin{array}{cc}
1 & a_{1/2} \\ 
-a_{1/2} & 1%
\end{array}%
\right) \oplus \left( 
\begin{array}{cc}
1 & a_{3/2} \\ 
-a_{3/2} & 1%
\end{array}%
\right) \oplus \ldots \oplus \left( 
\begin{array}{cc}
1 & a_{N-1/2} \\ 
-a_{N-1/2} & 1%
\end{array}%
\right) \right\} .  \label{example_K}
\end{equation}

That is, $K$ is a $2N$-by-$2N$ block-diagonal matrix that has blocks 
\begin{equation*}
\frac{1}{2}\left( 
\begin{array}{cc}
1 & a_{p} \\ 
-a_{p} & 1%
\end{array}%
\right)
\end{equation*}%
on its main diagonal.

It follows that it can be written as 
\begin{equation*}
K=\sum_{p}v_{p}v_{p}^{\ast },
\end{equation*}%
where 
\begin{equation*}
v_{p}^{\ast }=\frac{1}{\sqrt{2}}\left( 0,\ldots ,0,1,a_{p},0,\ldots
,0\right) ,
\end{equation*}%
and the entries $\ 1$ and $a_{p}$ are at places $2p$ and $2p+1,$
respectively.

Hence the kernel $\sigma _{N4}\left( \theta \right) $ has a quasi-real
diagonal form with eigenvalues $\lambda _{p}=1.$

Clearly the restriction of this random point field to an interval $I=\left(
a,b\right) $ is also a Pfaffian random point field with the kernel $\sigma
_{IN4}=\mathbf{1}_{I}\left( \theta \right) \sigma _{N4}\left( \theta -\theta
^{\prime }\right) \mathbf{1}_{I}\left( \theta ^{\prime }\right) $. Numerical
evaluations suggest that this restricted kernel is also quasi-real with
positive eigenvalues for arbitrary $N$ and $I.$ However, the proof of this
claim is elusive. The author proved it only for $N=2.$

\subsection{Gaussian Symplectic Ensemble}

The Gaussian symplectic ensembles of random matrices (GSE) consists of $N$%
-by-$N$ self-dual real quaternion matrices $H$ with the density given by the
formula 
\begin{equation*}
p\left( H\right) =c\exp \left[ -2\mathrm{Tr}\left( H^{2}\right) \right] .
\end{equation*}
The ensemble is called Gaussian because the entries of $H$ have components
that are real Gaussian variables with variance $1/4$.

The eigenvalues of a GSE matrix are real and have the density 
\begin{equation*}
c^{\prime }\prod\limits_{j<k}\left( x_{j}-x_{k}\right)
^{4}\prod\limits_{j=1}^{N}w\left( x_{j}\right) dx_{j},
\end{equation*}%
where $w\left( x\right) =\exp \left( -x^{2}\right) .$ Let $Q_{j}\left(
x\right) $ be polynomials of degree $j$ which are orthogonal with respect to
weight $w\left( x\right) .$ That is, : 
\begin{equation*}
\,\left\langle Q_{2j},Q_{2j+1}\right\rangle :=\int \left( Q_{2j}\left(
x\right) Q_{2j+1}^{\prime }\left( x\right) -Q_{2j}^{\prime }\left( x\right)
Q_{2j+1}\left( x\right) \right) w\left( x\right) dx=1
\end{equation*}%
$\left\langle Q_{2j+1},Q_{2j}\right\rangle =-1$ and $\left\langle
Q_{k},Q_{l}\right\rangle =0$ for all other choices of $k$ and $l.$

Define 
\begin{eqnarray*}
S_{N}\left( x,y\right) &=&\sqrt{w\left( x\right) w\left( y\right) }%
\sum_{k=0}^{2N-1}\left[ Q_{2k+1}^{\prime }\left( x\right) Q_{2k}\left(
y\right) -Q_{2k}^{\prime }\left( x\right) Q_{2k+1}\left( y\right) \right] ,
\\
I_{N}\left( x,y\right) &=&\sqrt{w\left( x\right) w\left( y\right) }%
\sum_{k=0}^{2N-1}\left[ Q_{2k+1}\left( x\right) Q_{2k}\left( y\right)
-Q_{2k}\left( x\right) Q_{2k+1}\left( y\right) \right] , \\
D_{N}\left( x,y\right) &=&\sqrt{w\left( x\right) w\left( y\right) }%
\sum_{k=0}^{2N-1}\left[ -Q_{2k+1}^{\prime }\left( x\right) Q_{2k}^{\prime
}\left( y\right) +Q_{2k}^{\prime }\left( x\right) Q_{2k+1}^{\prime }\left(
y\right) \right] ,
\end{eqnarray*}%
Define the quaternion kernel $K_{N}\left( x,y\right) $ by its matrix
representation: 
\begin{equation*}
\varphi \left( K_{N}\left( x,y\right) \right) =\left( 
\begin{array}{cc}
S_{N}\left( x,y\right) & D_{N}\left( x,y\right) \\ 
I_{N}\left( x,y\right) & S_{N}\left( y,x\right)%
\end{array}%
\right) .
\end{equation*}

The eigenvalues of GSE form a Pfaffian field with this kernel (see, for
example, Tracy and Widom \cite{tracy_widom98} or Chapter 5 in Mehta \cite%
{mehta04} for an explanation).

We can introduce the quaternion functions $\chi _{k}\left( x\right) $ as
follows: 
\begin{equation*}
\varphi \left( \chi _{k}\left( x\right) \right) =\sqrt{w\left( x\right) }%
\left( 
\begin{array}{cc}
Q_{2k}\left( x\right) & -Q_{2k}^{\prime }\left( x\right) \\ 
-Q_{2k+1}\left( x\right) & Q_{2k+1}^{\prime }\left( x\right)%
\end{array}%
\right) .
\end{equation*}%
The dual is 
\begin{equation*}
\varphi \left( \chi _{k}^{\ast }\left( x\right) \right) =\sqrt{w\left(
x\right) }\left( 
\begin{array}{cc}
Q_{2k+1}^{\prime }\left( x\right) & Q_{2k}^{\prime }\left( x\right) \\ 
Q_{2k+1}\left( x\right) & Q_{2k}\left( x\right)%
\end{array}%
\right) ,
\end{equation*}%
and we find that 
\begin{equation*}
K_{N}\left( x,y\right) =\sum_{k=0}^{N-1}\chi _{k}^{\ast }\left( x\right)
\chi _{k}\left( y\right) .
\end{equation*}%
Since $\chi _{k}$ are orthonormal, hence we find that the kernel $%
K_{N}\left( x,y\right) $ is finite rank and has a quasireal diagonal form
with all eigenvalues equal to $1.$

The kernel is positive since the determinants $\mathrm{Det}_{M}\left. \left(
K_{N}\left( x_{i},x_{j}\right) \right) \right\vert _{1\leq i,j\leq m}$ can
be interpreted as correlation functions for eigenvalues. Unfortunately, it
is not clear how to show the positivity without appeal to correlation
functions. For this reason, Proposition \ref{prop_projection_kernel} cannot
be used to give an independent proof that this kernel corresponds to a valid
random point field.

It is also an open question whether the restrictions of this kernel to
finite intervals have positive eigenvalues. This would be necessary to show
in order to establish the CLT by using Theorem \ref{theorem_CLT_general}.
(Note however that for the Gaussian symplectic ensemble, the CLT is already
known from the results in \cite{costin_lebowitz95}, which are based on the
relations between eigenvalues of Gaussian symplectic, orthogonal and unitary
ensembles.)

\bigskip

\appendix

\section{Quaternion matrices and determinants}

\label{section_determinants} The algebra of complex quaternions $\mathbb{Q}_{%
\mathbb{C}}$ is isomorphic to the algebra of two-by-two complex matrices $%
M_{2}\left( \mathbb{C}\right) ,$ with the correspondence defined by the rules%
\begin{equation*}
\mathbf{i}=\left( 
\begin{array}{cc}
0 & i \\ 
i & 0%
\end{array}%
\right) ,\text{ }\mathbf{j}=\left( 
\begin{array}{cc}
0 & -1 \\ 
1 & 0%
\end{array}%
\right) ,\text{ and }\mathbf{k}=\left( 
\begin{array}{cc}
i & 0 \\ 
0 & -i%
\end{array}%
\right) .
\end{equation*}%
In terms of $2$-by-$2$ matrices, if $q=\left( 
\begin{array}{cc}
a & b \\ 
c & d%
\end{array}%
\right) ,$ then its conjugate is $q^{\ast }=\left( 
\begin{array}{cc}
d & -b \\ 
-c & a%
\end{array}%
\right) .$

A number $\lambda $ is called an eigenvalue of a quaternion matrix $X$ if
for some non-zero quaternion vector $v$, we have $Xv=v\lambda .$ (These are
the right eigenvalues of $X,$ which are the most convenient in
applications.) It is easy to see that if $\lambda $ is an eigenvalue, then $%
q^{-1}\lambda q$ is also an eigenvalue for any quaternion $q.$ However, for
self-dual real quaternion matrices, all eigenvalues are real and it is
possible to show that every $n$-by-$n$ matrix $X$ of this type has exactly $%
n $ eigenvalues (counting with multiplicities); see Zhang \cite{zhang97}.

It is possible and useful to generalize the concept of determinant to
quaternion matrices. There are several sensible ways to do this and in this
paper we will only use the Moore-Dyson and Study determinants. Interested
reader can find details in a review paper by Aslaksen \cite{aslaksen96}.

If we replace each entry of a quaternion matrix $X$ by a corresponding $2$%
-by-$2$ complex matrix, then the matrix $X$ becomes represented by a $2n$-by-%
$2n$ complex matrix which we denote $\varphi \left( X\right) $. Then the
Study determinant of $X$ is defined as the usual determinant of $\varphi
\left( X\right) $. 
\begin{equation*}
\mathrm{Det}_{S}\left( X\right) :=\det \left( \varphi \left( X\right)
\right) .
\end{equation*}%
It can also be defined in a slightly different way for real quaternion
matrices. If $X=X_{1}+X_{2}\mathbf{i}+X_{3}\mathbf{j}+X_{4}\mathbf{k,}$
where $X_{1},$ $X_{2},$ $X_{3},$ and $X_{4}$ are real, then we define two
complex matrices $A=X_{1}+X_{2}i$ and $B=X_{3}+X_{4}i.$ Then, $\psi \left(
X\right) $ is defined as a $2n$-by-$2n$ complex matrix $\left( 
\begin{array}{cc}
A & B \\ 
-\overline{B} & \overline{A}%
\end{array}%
\right) ,$ where $\overline{A}$ is the conjugate of matrix $A,$ that is, $%
\left( \overline{A}\right) _{kl}=\overline{\left( A_{kl}\right) },$ and
similarly for $\overline{B}.$ The matrix $\psi \left( X\right) $ is called
the \emph{complex adjoint} of matrix $X.$ Then, $\mathrm{Det}_{S}\left(
X\right) =\det \left( \psi \left( X\right) \right) .$

The Study determinant is multiplicative: $\mathrm{Det}_{S}\left( AB\right) =%
\mathrm{Det}_{S}\left( A\right) \mathrm{Det}_{S}\left( B\right) $ for square
matrices $A$ and $B.$

The Moore-Dyson determinant of a self-dual real quaternion matrix $X$ can be
defined as the product of the right eigenvalues of the matrix. Remarkably,
this determinant can also be extended to all quaternion matrices by using a
variant of the Cayley combinatorial formula for the determinant. Namely, let 
$S_{n}$ be the group of permutations of the set $\left\{ 1,\ldots ,n\right\}
.$ Write every permutation $\sigma $ as a product of cycles: 
\begin{equation*}
\sigma =\left( n_{1}i_{2}\ldots i_{s}\right) \left( n_{2}j_{2}\ldots
j_{t}\right) \ldots \left( n_{r}k_{2}\ldots k_{l}\right) ,
\end{equation*}%
where $n_{i}$ are the largest elements of each cycle and $n_{1}>n_{2}>\ldots
>n_{r}.$ Then we can write 
\begin{equation*}
\mathrm{Det}_{M}\left( X\right) =\sum_{\sigma }\varepsilon \left( \sigma
\right) \left( X_{n_{1}i_{2}}X_{i_{2}i_{3}}\ldots X_{i_{s}n_{1}}\right)
\ldots \left( X_{n_{r}k_{2}}X_{k_{2}k_{3}}\ldots X_{k_{l}n_{r}}\right) ,
\end{equation*}%
where $\varepsilon \left( \sigma \right) =\left( -1\right) ^{n-r}$ is the
sign of the permutation $\sigma .$ (see \cite{moore22} and \cite{dyson70}).

Note that this definition allows one to calculate the quantity $\mathrm{Det}%
_{M}\left( X\right) $ for an arbitrary quaternion matrix. Dyson established
that this quantity is scalar for every self-dual matrix, that is, in this
case the $\mathbf{i},$ $\mathbf{j},$ and $\mathbf{k}$ components of the
determinant are zero.

The Moore-Dyson quaternion determinant of a self-dual quaternion matrix can
also be written as the Pfaffian of a related complex matrix. Let $J$ be a $%
2n $-by-$2n$ block-diagonal matrix with the blocks $\left( 
\begin{array}{cc}
0 & -1 \\ 
1 & 0%
\end{array}%
\right) $on the main diagonal. If $X$ is a self-dual quaternion matrix, then 
$-J\varphi \left( X\right) $ is antisymmetric (that is, $\left[ -J\varphi
\left( X\right) \right] ^{T}=J\varphi \left( X\right) $), and we can compute
the Pfaffian of this matrix. We have 
\begin{equation}
\mathrm{Det}_{M}\left( X\right) =\mathrm{Pf}\left( -J\varphi \left( X\right)
\right) .  \label{formula_determinant_pfaffian}
\end{equation}%
(see \cite{dyson70} and Proposition 6.1.5 on p. 238 in Forrester's book \cite%
{forrester10}).

In terms of the transformation $\psi ,$ this can be written as follows. Let $%
X$ be a real quaternion matrix and let $A$ and $B$ be defined as in the
definition of the complex adjoint. Let $\widetilde{J}=\psi \left( \mathbf{j}%
\right) =\left( 
\begin{array}{cc}
0 & -I \\ 
I & 0%
\end{array}%
\right) .$Then, $-\widetilde{J}\psi \left( X\right) =\left( 
\begin{array}{cc}
\overline{B} & \overline{A} \\ 
-A & B%
\end{array}%
\right) .$If the real quaternion matrix $X$ is self-dual, then $-\widetilde{J%
}\psi \left( X\right) $ is antisymmetric, and it can be shown that%
\begin{equation}
\mathrm{Det}_{M}\left( X\right) =-\mathrm{Pf}\left( -\widetilde{J}\psi
\left( X\right) \right) .  \label{det_pfaff}
\end{equation}

The Study and Moore-Dyson determinants are related by the following formula:

\begin{equation}
\mathrm{Det}_{S}\left( X\right) =\mathrm{Det}_{M}\left( X^{\ast }X\right) .
\label{identity_study_moore}
\end{equation}

(See formula (6.13) on page 239 in \cite{forrester10} or Corollary 5.1.3 on
p. 75 in \cite{mehta04}.)

\section{Proof of Theorem \protect\ref{theorem_cauchy_binet}}

Since the identity is algebraic, it is enough to show that it holds for
matrices with real quaternion entries. We will prove this by showing that
the corresponding result holds if we write the quaternion determinants in
terms of Pfaffians. Namely, let $C=X_{1}+X_{2}\mathbf{i}+X_{3}\mathbf{j}%
+X_{4}\mathbf{k}$ and define complex matrices $A=X_{1}+X_{2}i$ and $%
B=X_{3}+X_{4}i.$ Then, by using (\ref{det_pfaff}) we obtain that \ 
\begin{equation*}
\mathrm{Det}_{M}(C^{\ast }C)=-\mathrm{Pf}\left( 
\begin{array}{cc}
-B^{\ast }A+\left( B^{\ast }A\right) ^{t} & \left( A^{\ast }A\right)
^{t}+B^{\ast }B \\ 
-A^{\ast }A-\left( B^{\ast }B\right) ^{t} & A^{\ast }B-\left( A^{\ast
}B\right) ^{t}%
\end{array}%
\right) .
\end{equation*}

The blocks $-B^{\ast }A+\left( B^{\ast }A\right) ^{t}$ and $A^{\ast
}B-\left( A^{\ast }B\right) ^{t}$ are antisymmetric, and $\left( \left(
A^{\ast }A\right) ^{t}+B^{\ast }B\right) ^{t}=-A^{\ast }A-\left( B^{\ast
}B\right) ^{t},$ so the block matrix is antisymmetric as well.

What we need to prove is that 
\begin{equation*}
\mathrm{Pf}\left( 
\begin{array}{cc}
-B^{\ast }A+\left( B^{\ast }A\right) ^{t} & \left( A^{\ast }A\right)
^{t}+B^{\ast }B \\ 
-A^{\ast }A-\left( B^{\ast }B\right) ^{t} & A^{\ast }B-\left( A^{\ast
}B\right) ^{t}%
\end{array}%
\right) =\sum_{I}\mathrm{Pf}\left( 
\begin{array}{cc}
-B^{I\ast }A^{I}+\left( B^{I\ast }A^{I}\right) ^{t} & \left( A^{I\ast
}A^{I}\right) ^{t}+B^{I\ast }B^{I} \\ 
-A^{I\ast }A^{I}-\left( B^{I\ast }B^{I}\right) ^{t} & A^{I\ast }B^{I}-\left(
A^{I\ast }B^{I}\right) ^{t}%
\end{array}%
\right) ,
\end{equation*}%
where the summation is over all ordered $m$-tuples $I=\left( i_{1}<\ldots
<i_{m}\right) ,$ with $i_{k}\in \{1,\ldots ,n\},$ and $A^{I},$ $B^{I}$ are
the matrices that are obtained from matrices $A$ and $B,$ respectively, by
taking the rows with indices in $I.$

In order to prove this, we recall that if $R$ is a $2m$-by-$2m$
antisymmetric matrix, then the pfaffian of $R$ is defined as follows: 
\begin{equation}
\mathrm{Pf}\left( R\right) =\frac{1}{2^{m}m!}\sum_{\sigma \in
S_{2m}}sgn\left( \sigma \right) \prod_{i=1}^{m}R_{\sigma \left( 2i-1\right)
\sigma \left( 2i\right) }.  \label{def_pfaffian}
\end{equation}

Next, we note that if 
\begin{equation*}
R=\left( 
\begin{array}{cc}
-B^{\ast }A+\left( B^{\ast }A\right) ^{t} & \left( A^{\ast }A\right)
^{t}+B^{\ast }B \\ 
-A^{\ast }A-\left( B^{\ast }B\right) ^{t} & A^{\ast }B-\left( A^{\ast
}B\right) ^{t}%
\end{array}%
\right) ,
\end{equation*}%
then there is a formula for $R_{ij}$ in terms of elements of $A$ and $B.$
This formula depends on whether $i$ and $j$ are greater or less than $m.$
For example if $i$ and $j$ are both $\leq m,$ then 
\begin{equation*}
R_{ij}=\sum_{a=1}^{n}\left( -\overline{B}_{a,i}A_{a,j}+\overline{B}%
_{a,j}A_{a,i}\right) \equiv \sum_{a=1}^{n}\Psi _{a}\left( i,j\right) .
\end{equation*}

If we substitute this in formula (\ref{def_pfaffian}), and expand, then we
get 
\begin{equation}
\mathrm{Pf}\left( R\right) =\frac{1}{2^{m}m!}\sum_{\sigma \in
S_{2m}}sgn\left( \sigma \right) \sum_{\left( a,b,\ldots ,z\right) }\Psi
_{a}\left( \sigma \left( 1\right) ,\sigma \left( 2\right) \right) \ldots
\Psi _{z}\left( \sigma \left( 2m-1\right) ,\sigma \left( 2m\right) \right) ,
\label{expansion_R}
\end{equation}%
where the summation is over all $m$-tuples $\left( a,b,\ldots ,z\right) $
with each letter taking a value in $\left\{ 1,\ldots ,n\right\} .$

A similar formula holds for the pfaffian of $R^{I},$ where 
\begin{equation*}
R^{I}=\left( 
\begin{array}{cc}
-B^{I\ast }A^{I}+\left( B^{I\ast }A^{I}\right) ^{t} & \left( A^{I\ast
}A^{I}\right) ^{t}+B^{I\ast }B^{I} \\ 
-A^{I\ast }A^{I}-\left( B^{I\ast }B^{I}\right) ^{t} & A^{I\ast }B^{I}-\left(
A^{I\ast }B^{I}\right) ^{t}%
\end{array}%
\right) .
\end{equation*}%
Namely, 
\begin{equation}
\mathrm{Pf}\left( R^{I}\right) =\frac{1}{2^{m}m!}\sum_{\sigma \in
S_{2m}}sgn\left( \sigma \right) \sum_{\left( a,b,\ldots ,z\right) \in
I^{m}}\Psi _{a}\left( \sigma \left( 1\right) ,\sigma \left( 2\right) \right)
\ldots \Psi _{z}\left( \sigma \left( 2m-1\right) ,\sigma \left( 2m\right)
\right) .  \label{expansion_RI}
\end{equation}%
The difference with the previous formula is that the elements of $m$-tuples $%
\left( a,b,\ldots ,z\right) $ are now restricted to take values among
indices in $m$-tuple $I.$

Let us for shortness write $\Psi _{a,\ldots ,z}\left( \sigma \right) $ for
the product $\Psi _{a}\left( \sigma \left( 1\right) ,\sigma \left( 2\right)
\right) \ldots \Psi _{z}\left( \sigma \left( 2m-1\right) ,\sigma \left(
2m\right) \right) .$

If all elements of $\left( a,b,\ldots ,z\right) $ are different, then the
term $\Psi _{a,\ldots ,z}\left( \sigma \right) $ occurs once in expansion (%
\ref{expansion_R}) and once in the sum of expansions (\ref{expansion_RI}), 
\begin{equation*}
\sum_{I}\mathrm{Pf}\left( R^{I}\right) .
\end{equation*}%
(In this sum, it occurs in the expansion of that $\mathrm{Pf}\left(
R^{I}\right) ,$ for which $I$ is the ordered version of the $m$-tuple $%
\left( a,b,\ldots ,z\right) .$)

If some of the elements of $\left( a,b,\ldots ,z\right) $ coincide, then the
situation is different. The term $\Psi _{a,\ldots ,z}\left( \sigma \right) $
occurs once in expansion (\ref{expansion_R}) but it can occur more than once
in the sum 
\begin{equation*}
\sum_{I}\mathrm{Pf}\left( R^{I}\right) .
\end{equation*}%
For example, if all elements of the $m$-tuples are the same, $a=b=\ldots =z,$
then this term will appear in the expansion of each $\mathrm{Pf}\left(
R^{I}\right) ,$ whose index $I$ contain $a.$

Clearly, in order to prove that $\mathrm{Pf}\left( R\right) =\sum_{I}\mathrm{%
Pf}\left( R^{I}\right) ,$ it is enough to prove that the sum of all these
terms is zero. That is, it is enough to show that for a fixed $m$-tuple $%
\left( a,b,\ldots ,z\right) $ with at least two elements that are equal, the
sum 
\begin{equation*}
\sum_{\sigma \in S_{2m}}sgn\left( \sigma \right) \Psi _{a,\ldots ,z}\left(
\sigma \right)
\end{equation*}%
is zero.

Without loss of generality we can assume that $a=b.$ We have to consider
several cases of $\sigma ,$ which are summarized in the following table

\begin{tabular}{ccccc}
& $\sigma \left( 1\right) $ & $\sigma \left( 2\right) $ & $\sigma \left(
3\right) $ & $\sigma \left( 4\right) $ \\ \hline
$S_{2m}\left[ 1\right] $ &  &  &  &  \\ \hline
$S_{2m}\left[ 2\right] $ &  &  & * & * \\ \hline
$S_{2m}\left[ 3\right] $ &  & * &  & * \\ \hline
$S_{2m}\left[ 4\right] $ &  & * & * &  \\ \hline
$S_{2m}\left[ 5\right] $ &  &  &  & * \\ \hline
$S_{2m}\left[ 6\right] $ &  & * &  &  \\ \hline
$S_{2m}\left[ 7\right] $ &  &  & * &  \\ \hline
$S_{2m}\left[ 8\right] $ &  & * & * & *%
\end{tabular}

The star means that the corresponding $\sigma \left( i\right) $ is greater
than $m.$ For example, $S_{2m}\left[ 1\right] $ denote the set of all
permutations from $S_{2m}$ that satisfy the condition that all of $\sigma
\left( 1\right) ,\sigma \left( 2\right) ,\sigma \left( 3\right) ,\sigma
\left( 4\right) $ are smaller than or equal to $m.$ For permutations in this
set, 
\begin{equation*}
\Psi _{a}\left( \sigma \left( 1\right) ,\sigma \left( 2\right) \right) \Psi
_{a}\left( \sigma \left( 3\right) ,\sigma \left( 4\right) \right) =\left( -%
\overline{B}_{a,i}A_{a,j}+\overline{B}_{a,j}A_{a,i}\right) \left( -\overline{%
B}_{a,k}A_{a,l}+\overline{B}_{a,l}A_{a,k}\right) .
\end{equation*}

$S_{2m}\left[ 2\right] $ denote the set of all permutations from $S_{2m}$
that satisfy the condition that $\sigma \left( 1\right) ,\sigma \left(
2\right) $ are smaller than or equal to $m$ and $\sigma \left( 3\right)
,\sigma \left( 4\right) $ are greater than $m,$ and so on.

Let us define $\tau _{1}\left[ \sigma \right] ,$ as a permutation that
coincides with $\sigma $ on all indices except $2$ and $4,$ for which it is
defined by equalities $\tau _{1}\left[ \sigma \right] \left( 2\right)
=\sigma \left( 4\right) ,$ and $\tau _{1}\left[ \sigma \right] \left(
4\right) =\sigma \left( 2\right) .$ Similarly, $\tau _{2}\left[ \sigma %
\right] $ is defined as a permutation which acts on everything except $2$
and $3$ as $\sigma ,$ and on these indices it is defined by $\tau _{2}\left[
\sigma \right] \left( 2\right) =\sigma \left( 3\right) ,$ $\tau _{2}\left[
\sigma \right] \left( 3\right) =\sigma \left( 2\right) .$ Finally we define $%
\tau _{3}\left[ \sigma \right] $ as a permutation that coincides with $%
\sigma $ on all indices except $2,3,$ and $4,$ where it is defined by the
rules: $\tau _{3}\left[ \sigma \right] \left( 2\right) =\sigma \left(
4\right) ,$ $\tau _{3}\left[ \sigma \right] \left( 3\right) =\sigma \left(
2\right) ,$ $\tau _{3}\left[ \sigma \right] \left( 4\right) =\sigma \left(
3\right) .$ Note that $sgn\left( \tau _{1}\left[ \sigma \right] \right)
=sgn\left( \tau _{2}\left[ \sigma \right] \right) =-sgn\left( \sigma \right)
,$ and $sgn\left( \tau _{3}\left[ \sigma \right] \right) =sgn\left( \sigma
\right) .$ Observe that for an arbitrary function $f,$ 
\begin{equation*}
\sum_{\sigma \in S_{2m}\left[ 1\right] }\left( f\left( \sigma \right)
+f\left( \tau _{1}\left[ \sigma \right] \right) +f\left( \tau _{2}\left[
\sigma \right] \right) \right) =3\sum_{\sigma \in S_{2m}\left[ 1\right]
}f\left( \sigma \right) .
\end{equation*}%
It is easy to check that the identity 
\begin{eqnarray*}
0 &=&\Psi _{a}\left( \sigma \left( 1\right) ,\sigma \left( 2\right) \right)
\Psi _{a}\left( \sigma \left( 3\right) ,\sigma \left( 4\right) \right) -\Psi
_{a}\left( \tau _{1}\left[ \sigma \right] \left( 1\right) ,\tau _{1}\left[
\sigma \right] \left( 2\right) \right) \Psi _{a}\left( \tau _{1}\left[
\sigma \right] \left( 3\right) ,\tau _{1}\left[ \sigma \right] \left(
4\right) \right) \\
&&-\Psi _{a}\left( \tau _{2}\left[ \sigma \right] \left( 1\right) ,\tau _{2}%
\left[ \sigma \right] \left( 2\right) \right) \Psi _{a}\left( \tau _{2}\left[
\sigma \right] \left( 3\right) ,\tau _{2}\left[ \sigma \right] \left(
4\right) \right)
\end{eqnarray*}%
holds for permutations in $S_{2m}\left[ 1\right] $, and this implies that 
\begin{equation}
\sum_{\sigma \in S_{2m}\left[ 1\right] }sgn\left( \sigma \right) \Psi
_{a,\ldots ,z}\left( \sigma \right) =0.  \label{sum_1}
\end{equation}

Next, 
\begin{equation*}
\sum_{\sigma \in S_{2m}\left[ 2\right] }\left( f\left( \sigma \right)
+f\left( \tau _{1}\left[ \sigma \right] \right) +f\left( \tau _{3}\left[
\sigma \right] \right) \right) =\sum_{\sigma \in S_{2m}\left[ 2\right]
}f\left( \sigma \right) +2\sum_{\sigma \in S_{2m}\left[ 3\right] }f\left(
\sigma \right) .
\end{equation*}%
and it is easy to check the identity 
\begin{eqnarray}
0 &=&\Psi _{a}\left( \sigma \left( 1\right) ,\sigma \left( 2\right) \right)
\Psi _{a}\left( \sigma \left( 3\right) ,\sigma \left( 4\right) \right) -\Psi
_{a}\left( \tau _{1}\left[ \sigma \right] \left( 1\right) ,\tau _{1}\left[
\sigma \right] \left( 2\right) \right) \Psi _{a}\left( \tau _{1}\left[
\sigma \right] \left( 3\right) ,\tau _{1}\left[ \sigma \right] \left(
4\right) \right)  \notag \\
&&+\Psi _{a}\left( \tau _{3}\left[ \sigma \right] \left( 1\right) ,\tau _{3}%
\left[ \sigma \right] \left( 2\right) \right) \Psi _{a}\left( \tau _{3}\left[
\sigma \right] \left( 3\right) ,\tau _{3}\left[ \sigma \right] \left(
4\right) \right) .  \label{identity_case3}
\end{eqnarray}%
for $\sigma \in S_{2m}\left[ 2\right] .$ Hence, identity (\ref%
{identity_case3}) implies that%
\begin{equation}
\left( \sum_{\sigma \in S_{2m}\left[ 2\right] }+2\sum_{\sigma \in S_{2m}%
\left[ 3\right] }\right) sgn\left( \sigma \right) \Psi _{a,\ldots ,z}\left(
\sigma \right) =0  \label{identity_case3_B}
\end{equation}

The other cases are similar. We use

\begin{equation*}
\sum_{\sigma \in S_{2m}\left[ 4\right] }\left( f\left( \sigma \right)
+f\left( \tau _{1}\left[ \sigma \right] \right) +f\left( \tau _{3}\left[
\sigma \right] \right) \right) =2\sum_{\sigma \in S_{2m}\left[ 4\right]
}f\left( \sigma \right) +\sum_{\sigma \in S_{2m}\left[ 2\right] }f\left(
\sigma \right)
\end{equation*}%
in order to conclude that 
\begin{equation}
\left( 2\sum_{\sigma \in S_{2m}\left[ 4\right] }+\sum_{\sigma \in S_{2m}%
\left[ 2\right] }\right) sgn\left( \sigma \right) \Psi _{a,\ldots ,z}\left(
\sigma \right) =0  \label{identity_case4_B}
\end{equation}

By adding (\ref{identity_case3_B}) and (\ref{identity_case4_B}), we obtain%
\begin{equation}
\left( \sum_{\sigma \in S_{2m}\left[ 2\right] }+\sum_{\sigma \in S_{2m}\left[
3\right] }+\sum_{\sigma \in S_{2m}\left[ 4\right] }\right) sgn\left( \sigma
\right) \Psi _{a,\ldots ,z}\left( \sigma \right) =0.  \label{sum_234}
\end{equation}

Next, the identity 
\begin{equation*}
\sum_{\sigma \in S_{2m}\left[ 5\right] }\left( f\left( \sigma \right)
+f\left( \tau _{1}\left[ \sigma \right] \right) +f\left( \tau _{3}\left[
\sigma \right] \right) \right) =2\sum_{\sigma \in S_{2m}\left[ 5\right]
}f\left( \sigma \right) +\sum_{\sigma \in S_{2m}\left[ 6\right] }f\left(
\sigma \right)
\end{equation*}%
implies that 
\begin{equation}
\left( 2\sum_{\sigma \in S_{2m}\left[ 5\right] }+\sum_{\sigma \in S_{2m}%
\left[ 6\right] }\right) sgn\left( \sigma \right) \Psi _{a,\ldots ,z}\left(
\sigma \right) =0,  \label{identity_case4_C}
\end{equation}%
and the identity 
\begin{equation*}
\sum_{\sigma \in S_{2m}\left[ 6\right] }\left( f\left( \sigma \right)
+f\left( \tau _{1}\left[ \sigma \right] \right) +f\left( \tau _{3}\left[
\sigma \right] \right) \right) =\sum_{\sigma \in S_{2m}\left[ 6\right]
}f\left( \sigma \right) +2\sum_{\sigma \in S_{2m}\left[ 7\right] }f\left(
\sigma \right)
\end{equation*}%
implies that 
\begin{equation}
\left( \sum_{\sigma \in S_{2m}\left[ 6\right] }+2\sum_{\sigma \in S_{2m}%
\left[ 7\right] }\right) sgn\left( \sigma \right) \Psi _{a,\ldots ,z}\left(
\sigma \right) =0.  \label{identity_case4_D}
\end{equation}

By adding (\ref{identity_case4_C}) and (\ref{identity_case4_D}), we obtain: 
\begin{equation}
\left( \sum_{\sigma \in S_{2m}\left[ 5\right] }+\sum_{\sigma \in S_{2m}\left[
6\right] }+\sum_{\sigma \in S_{2m}\left[ 7\right] }\right) sgn\left( \sigma
\right) \Psi _{a,\ldots ,z}\left( \sigma \right) =0.  \label{sum_567}
\end{equation}%
Finally, 
\begin{equation*}
\sum_{\sigma \in S_{2m}\left[ 8\right] }\left( f\left( \sigma \right)
+f\left( \tau _{1}\left[ \sigma \right] \right) +f\left( \tau _{3}\left[
\sigma \right] \right) \right) =3\sum_{\sigma \in S_{2m}\left[ 8\right]
}f\left( \sigma \right) .
\end{equation*}%
Identity (\ref{identity_case3}) still holds in this case, and therefore, 
\begin{equation}
\sum_{\sigma \in S_{2m}\left[ 8\right] }sgn\left( \sigma \right) \Psi
_{a,\ldots ,z}\left( \sigma \right) =0.  \label{sum_8}
\end{equation}

Clearly, the sets $S_{2m}\left[ k\right] ,$ are disjoint and their union
over $k=1,\ldots ,8$ consists of all permutations from $S_{2m}$ for which $%
\sigma \left( 1\right) \leq m.$ Hence, identities (\ref{sum_1}), (\ref%
{sum_234}), (\ref{sum_567}), and (\ref{sum_8}) imply that 
\begin{equation*}
\sum_{\substack{ \sigma \in S_{2m}  \\ \sigma \left( 1\right) \leq m}}%
sgn\left( \sigma \right) \Psi _{a,\ldots ,z}\left( \sigma \right) =0.
\end{equation*}

The case $\sigma \left( 1\right) >m$ can be handled similarly, and we obtain%
\begin{equation*}
\sum_{\sigma \in S_{2m}}sgn\left( \sigma \right) \Psi _{a,\ldots ,z}\left(
\sigma \right) =0.
\end{equation*}%
This holds provided the first two indices in $\left( a,\ldots ,z\right) $
are equal. It is clear that this identity also holds if any two indices in $%
\left( a,\ldots ,z\right) $ are equal.

As observed earlier, this implies that we can remove all terms $\Psi
_{a,\ldots ,z}\left( \sigma \right) ,$ for which two indices in $\left(
a,\ldots ,z\right) $ coincide, from the expansions of both $\mathrm{Pf}%
\left( R\right) $ and the sum $\sum_{I}\mathrm{Pf}\left( R^{I}\right) .$
(See expansions (\ref{expansion_R}) and (\ref{expansion_RI}).) Then it is
clear that the remaining terms in the expansions are the same. Hence 
\begin{equation*}
\mathrm{Pf}\left( R\right) =\sum_{I}\mathrm{Pf}\left( R^{I}\right) ,
\end{equation*}%
and this implies the statement of the theorem. $\square $

\bibliographystyle{plain}
\bibliography{comtest}

\end{document}